# CARACTERIZACIÓN SEMÁNTICO-DEDUCTIVA DE LA LÓGICA DOBLE LD Y LOS GRÁFICOS EXISTENCIALES GAMMA-LD
# SEMANTIC AND DEDUCTIVE CHARACTERIZATION OF DOUBLE LOGIC LD AND EXISTENTIAL GRAPHS GAMMA-LD

Manuel Sierra-Aristizábal[1]
orcid.org/0000-0001-8157-2577
arXiv.2310.02035

**Resumen.** En este trabajo se presenta el sistema deductivo Lógica Proposicional Doble, LD, junto con la semántica de mundos posibles que lo caracterizan. LD incluye un operador de afirmación alterna y uno de negación alterna y no vale el principio del tercero excluido cuando uno de los disyuntos es uno de los operadores alternos. El sistema LD incluye como teoremas, los teoremas de la Lógica Proposicional Clásica LC. En LD se recuperan los conectivos intuicionistas a partir de los clásicos y de los alternos, lo cual implica que, el sistema LD incluye como teoremas, los teoremas de la Lógica Proposicional Intuicionista LI. Se prueban, de manera rigurosa y detallada, los teoremas de consistencia, validez y completitud. Se ilustra la capacidad que tiene LD para solucionar una versión de la paradoja del mentiroso, donde LC y LI fracasan. Finalmente, se determina de manera precisa, la ubicación del sistema LD, en relación con los sistemas que permitieron intuir su estructura, tales como, la Lógica Básica para la Verdad Aristotélica LBVA, la Lógica de las Tautologías LT y la Lógica Básica para la Verdad y la Falsedad LBVF. LD se encuentra caracterizada mediante los Gráficos Existenciales Gamma-LD, en el estilo de Peirce, por lo que la semántica de mundos posibles presentada en este trabajo también caracteriza a Gamma-LD.

**Palabras clave**: afirmación alterna, falsedad, mundos posibles, negación alterna, verdad.

**Abstract.** This work presents the deductive system Double Propositional Logic, LD, along with the semantics of possible worlds that characterize it. LD includes an alternate affirmation operator and an alternate negation operator and is not valid for the principle of the excluded third when one of the disjunctions is one of the alternate operators. The LD system includes as theorems, the theorems of the Classic Propositional Logic LC. In LD, the intuitionist connectives are recovered from the classics and the alternates, which implies that, the LD system includes as theorems, the theorems of the Propositional Intuitionist Logic LI. It is rigorously and detailed to prove that theorems of consistency, validity and completeness. It illustrates LD's ability to fix a version of the liar's paradox, where LC and LI fail. Finally, the location of the LD system is precisely determined in relation to the systems that made it possible to intuit the structure of LD, such as, Basic Logic for the Aristotelian Truth LBVA, the Logic of Tautologies LT and the Basic Logic for Truth and Falsehood LBVF. LD is characterized by Gamma-LD Existential Graphs, in the style of Peirce, so the semantics of possible worlds presented in this work also characterizes Gamma-LD.

## 1 PRESENTACIÓN

En 2002 se presenta, como una extensión propia de la *Lógica Proposicional Clásica* LC, el sistema *Lógica Básica Paraconsistente y Paracompleta* LB (ver Sierra-Aristizábal (2010)), este sistema incluye operadores de negación y afirmación alternas, de tal manera que se puede controlar la generación de fragmentos de lógicas paraconsistentes y/o paracompletas. En 2007 se presenta, como una extensión propia del primero, el sistema *Lógica Básica para la Verdad Aristotélica* LBVA Sierra-Aristizábal (2007a), y como una extensión propia del segundo, el sistema *Lógica Básica para la Verdad y la Falsedad* LBVF Sierra-Aristizábal (2007b).

Se presenta en Sierra-Aristizábal (2012) el sistema *Lógica de las Tautologías* LT, el cual es una extensión propia de LBVA, pero es independiente de LBVF. LT es caracterizado por una semántica de mundos posibles, en contraste con sistemas previos que son caracterizados por semánticas de valuaciones. Los últimos dos sistemas son bastante fuertes y tienen la característica de que no vale ley del tercero excluido, cuando uno de los disyuntos es una afirmación alterna o una negación alterna, lo cual genera algunos resultados típicos de la *Lógica*

---
[1] *msierra@eafit.edu.co*

*Proposicional Intuicionista* LI (ver Van Dale (2013)). Una pregunta natural en este punto es ¿Cómo está conectado el sistema LI, o una generalización del mismo, con los 4 sistemas mencionados?

En este trabajo se presenta el sistema deductivo *Lógica Doble* LD, como una extensión generalizada de la lógica proposicional clásica, la cual incluye fórmulas atómicas alternas y operadores alternos de afirmación y de negación. A partir de estos operadores, en LD se definen operadores alternos para los conectivos binarios clásicos, resultando que LD también es una extensión generalizada de LI. LD resulta ser independiente de LBVF, y además es un sistema intermedio entre LBVA y LT. En Sierra (2022) se presenta la caracterización de LD mediante los Gráficos Existenciales Gamma-LD, en el estilo de Peirce (1992), por lo que la semántica de mundos posibles presentada en este trabajo también caracteriza a Gamma-LD.

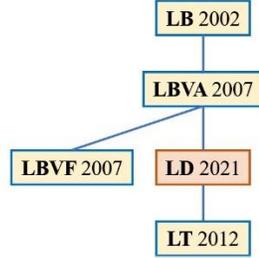

Al ser interpretado el operador de afirmación alterna como *verdad Aristotélica*, y el operador de negación alterna como *Falsedad Aristotélica*, en LD se tienen como consecuencia, las definiciones de verdad y falsedad dadas en Aristóteles (1998): "Decir de lo que es que es, y de lo que no es que no es, es lo verdadero; decir de lo que es que no es, y de lo que no es que es, es lo falso"; o previamente en Platón (1983): "El discurso, que dice las cosas como son, es verdadero; y el que las dice como no son, es falso". Con esta lectura en LD, así como en LBVF, se da solución a una de las versiones de la *paradoja del mentiroso* (ver Bochenski (1976)).

Respecto a la semántica de mundos posibles con la que se caracteriza LD, las pruebas de validez, consistencia y completitud son presentadas de manera detallada, resultando que LD es una extensión propia de la Lógica Modal S4 presentada en Lewis and Langford (1932). Finalmente se muestra, que LC y LI son casos particulares de LD.

## 2 LENGUAJE

Se parte de un conjunto infinito enumerable de **Fórmulas Atómicas Clásicas FAC**, las cuales se denotan $a, b, c, d, \ldots, x, y, z, \ldots$. A partir de FAC se obtiene el conjunto de **Fórmulas Atómicas Alternas FAA**, las cuales se denotan $\underline{a}, \underline{b}, \underline{c}, \underline{d}, \ldots, \underline{x}, \underline{y}, \underline{z}, \ldots$. El conjunto de **Fórmulas Atómicas FAT**, se obtiene de la unión de los dos anteriores FAT =FAC∩FAA.

El conjunto **FC** de **Fórmulas Clásicas** se define de la siguiente manera:
1. a∈FAC ⟹ a∈FC.
2. X∈FC ⟹ ∼X∈FC
3. X,Y∈FC ⟹ X⊃Y, X∪Y, X•Y, X≡Y∈FC.
4. Sólo 1 a 3 determinan FC.

El conjunto **FI** de **Fórmulas Intuicionistas** se define de la siguiente manera:
1. $\underline{a}$∈FAA ⟹ $\underline{a}$∈FI
2. X∈FI ⟹ ¬X∈FI.
3. X,Y∈FI ⟹ X→Y, X∨Y, X∧Y, X↔Y∈FI.
4. Sólo 1 a 3 determinan FI.

El conjunto **FOR** de **Fórmulas de LD** se define de la siguiente manera:
1. X∈FC ⟹ X∈FOR.
2. X∈FI ⟹ X∈FOR.
3. X∈FOR ⟹ ∼X∈FOR.
4. X∈FOR ⟹ ¬X∈FOR.
5. X,Y∈FOR ⟹ X⊃Y, X∪Y, X•Y, X≡Y∈FOR.
6. X,Y∈FOR ⟹ X→Y, X∨Y, X∧Y, X↔Y∈FOR.
7. Sólo 1 a 6 determinan FOR.

El conjunto **FA** de **Fórmulas alternas** se define de la siguiente manera:
$X \in FA \Leftrightarrow X = \neg Y$ con $Y \in FOR$ o $X = \underline{a}$ donde $\underline{a} \in FAA$

**Notación.** Si $X \in FOR$ entonces $\underline{X}$ indica que $X \in FA$.

Definición 1. (*Operadores de verdad*). En **LD**, a partir de la *Falsedad alterna* $\neg X$, se definen los conectivos (donde $X, Y \in FOR$):

*Verdad alterna*: $+X = \neg \sim X$. *Disyunción alterna*: $X \vee Y = +(X \cup Y)$.
*Conjunción alterna*: $X \wedge Y = +(X \bullet Y)$. *Implicación alterna*: $X \rightarrow Y = +(X \supset Y)$.
*Equivalencia alterna*: $X \leftrightarrow Y = +(X \equiv Y)$. *Verdad alterna refutable*: $-X = \sim +X$.
*Falsedad alterna refutable*: $\otimes X = \sim \neg X$. *Buena fundamentación*: $*X = (\neg X \cup +X)$.

Como consecuencias se tienen las siguientes equivalencias entre los llamados **operadores de verdad**:
a. $+X \equiv \neg \sim X \equiv \sim \otimes \sim X \equiv \sim -X$. b. $\neg X \equiv +\sim X \equiv \sim \otimes X \equiv \sim -\sim X$.
c. $\otimes X \equiv \sim + \sim X \equiv \sim \neg X \equiv - \sim X$. d. $-X \equiv \sim +X \equiv \sim \neg \sim X \equiv \otimes \sim X$.
e. $*X \equiv +X \cup \neg X \equiv -X \supset \neg X \equiv \otimes X \supset +X$.

## 3 SISTEMA DEDUCTIVO

El *sistema deductivo para* **LD** consta de los siguientes axiomas (donde $X, Y, Z \in FOR$, $\underline{X} \in FA$):
Ax1.1 $X \supset (Y \supset X)$  Ax1.2 $(X \supset (Y \supset Z)) \supset ((X \supset Y) \supset (X \supset Z))$
Ax1.3 $X \supset (X \cup Y)$  Ax1.4 $Y \supset (X \cup Y)$
Ax1.5 $(X \supset Z) \supset ((Y \supset Z) \supset ((X \cup Y) \supset Z))$  Ax1.6 $(X \bullet Y) \supset X$
Ax1.7 $(X \bullet Y) \supset Y$  Ax1.8 $(X \supset Y) \supset ((X \supset Z) \supset (X \supset (Y \bullet Z)))$
Ax1.9 $X \supset (\sim X \supset Y)$  Ax1.10 $X \cup \sim X$
Ax1.11 $(X \equiv Y) \supset (X \supset Y)$  Ax1.12 $(X \equiv Y) \supset (Y \supset X)$
Ax1.13 $(X \supset Y) \supset [(Y \supset X) \supset (X \equiv Y)]$
Ax2.1 AxMP+: $+(X \supset Y) \supset (+X \supset +Y)$  Ax2.2 AxR: $\neg X \supset \sim X$
Ax2.3 AxT: $\underline{X} \supset +\underline{X}$  Ax2.4 Ax+: $X \in \{Ax1.1, \ldots, Ax2.3\} \Rightarrow +X$ es un axioma.
Como única *regla de inferencia* se tiene el *Modus Ponens* MP: de X y $X \supset Y$ se infiere Y.

Definición 2. (*Tipos de axiomas*). Axiomas clásicos AxCP = {Ax1.1, …, Ax1.13}.
Axiomas alternos AxAL = {Ax2.1, …, Ax2.4}. Axiomas dobles AxLD = AxCP $\cap$ AxAL.

Utilizando los operadores de verdad se infiere que:
Ax2.2 AxR: $\neg X \supset \sim X$, puede ser presentado como AxR: $+X \supset X$, o como AxR: $X \supset \otimes X$.
Ax2.3 AxT: $\underline{X} \supset +\underline{X}$, puede ser presentado como AxT: $+X \supset ++X$, o como AxT: $\neg X \supset +\neg X$,
o como AxT: $-\neg X \supset \otimes X$.

Definición 3. (*Teorema*). Para $X \in FOR$. Se dice que X es un *teorema* ($X \in$ **TEO**) si y solamente si existe una *demostración de X*, es decir, X es la última fórmula de una sucesión finita de fórmulas, tales que cada una de ellas es un axioma o se infiere de dos fórmulas anteriores utilizando la regla de inferencia MP. El número de elementos de la sucesión se llama la *longitud de la prueba.*

Proposición 1 (*Construcción de verdades Alternas*). Para $X \in FOR$.
$X \in TEO \Rightarrow +X \in TEO$.
Prueba. Supóngase que X es un teorema, se probará que $+X \in TEO$, haciendo inducción sobre la longitud de la demostración de X.

**PB**. Paso base. La longitud de la demostración de X es 1, es decir X es un axioma, pero si X es un axioma entonces, por Ax+, +X es axioma y por lo tanto, +X∈TEO.

**PI**. Paso de inducción. Como hipótesis inductiva **HI**, se tiene que si la longitud de la demostración de Y es menor que L entonces +Y es un teorema. Supóngase que la demostración de X tiene longitud L>1. Se tiene entonces que X es un axioma o X es consecuencia de pasos anteriores utilizando la regla de inferencia MP. En el primer caso se procede como en el paso base. En el segundo caso se tienen, para alguna fórmula Z, demostraciones de Z→X y de Z, ambas de longitud menor que L. Aplicando **HI** resulta que +(Z→X),+Z∈TEO. Por AxMP+ se tiene +(Z→X)→(+Z→+X), aplicando dos veces la regla MP resulta que +X∈TEO. Por lo que, según el principio de inducción matemática, se ha probado la proposición. □

En lo que sigue se utilizaran los siguientes resultados del cálculo proposicional clásico CP, y como en LD se tienen los correspondientes axiomas de CP y la regla de inferencia MP entonces, estos resultados valen en LD (para detalles de las pruebas en CP ver Caicedo (1990) y Hamilton (1981)). Sean X,Y,Z,W∈FOR.

I•c. *Introducción*: ((Z⊃X)•(Z⊃Y))⊃[Z⊃(X•Y)).
I•. *Introducción de* •: X⊃(Y⊃(X•Y)).
E•c. *Eliminación*: (Z⊃(X•Y))⊃((Z⊃X)•(Z⊃Y)).
E•. *Eliminación de* •: (X•Y)⊃X, (X•Y)⊃Y.
SH. *Silogismo hipotético*: ((X⊃Y)•(Y⊃Z))⊃(X⊃Z).
Exp. *Exportación*: (X⊃(Y⊃Z))≡( (X•Y)⊃Z).
EQ. *Equivalencia*: (X≡Y)≡((X⊃Y)•(Y⊃X)).
DN. *Doble negación*: X≡~~X.
DI. *Demostración indirecta*: X⊃(Y•~Y) ⇒ ~X.
N∪. *Negación de* ∪: ~(X∪Y)≡(~X•~Y).
N•. *Negación de* •: ~(X•Y)≡(~X∪~Y).
N⊃. *Negación de* ⊃: ~(X⊃Y)≡(X•~Y).
Tras. *Transposición*: (X⊃Y)≡(~Y→~X).
Imp. *Implicación*: (X⊃Y)≡(~X∪Y).
Id. *Principio de identidad*: X⊃X.
PR. *Principio de retorsión*: de (~X⊃X)≡X.
TD. *Teorema de deducción*: (X⇒Y) ⇒ X⊃Y).
EQ. *Equivalencia*: (X≡Y)≡((X•Y)∪(~X•~Y)).
SustEQ. F(X),X,Y∈FC y X≡Y ⇒ F(X)⊃F(Y).

**Proposición 2** (*Conjunción de verdades Alternas*). Para X,Y∈FOR.

a. +(X•Y)≡(+X•+Y)∈TEO.
b. ¬(X∪Y)≡(¬X•¬Y)∈TEO.
c. (¬X∪¬Y)⊃¬(X•Y)∈TEO.

Prueba. Parte a. Por Axn.6, y Axn.7 se tienen (X•Y)⊃X y (X•Y)⊃Y, y por Ax+ resultan +((X•Y)⊃X) y +((X•Y)⊃Y). Utilizando el AxMP+ y MP se infieren +(X•Y)⊃+X y +(X•Y)⊃+Y. Utilizando I•c se infiere +(X•Y)⊃(+X•+Y).

Por otro lado, de la regla I• se tiene X⊃(Y⊃(X•Y)). Por la proposición 1 se infiere +(X⊃(Y⊃(X•Y))), y utilizando AxMP+ y MP resulta +X⊃+(Y⊃(X•Y)), como además por AxMP+ se tiene +(Y⊃(X•Y))⊃ (+Y⊃+(X•Y)), entonces por SH se obtiene +X⊃(+Y⊃ +(X•Y)). Utilizando Exp se infiere (+X •+Y)⊃+(X•Y). Como ya se probó la recíproca, entonces por I• y EQ resulta que +(X•Y)≡(+X•+Y)∈TEO.

Parte b. Por N∪ se tienen ~(X∪Y)⊃(~X•~Y) y (~X•~Y)⊃~(X∪Y), por la proposición 1, AxMP+ y EQ se infiere que +~(X∪Y)≡+(~X•~Y), por la parte a resulta +~(X∪Y)≡(+~X•+~Y), lo cual por falsedad alterna significa que ¬(X∪Y)≡(¬X•¬Y)∈TEO. □

Parte c. Por Ax1.6 y Ax1.7 se tienen X•Y⊃X∈TEO y X•Y⊃Y∈TEO, por Tras resultan ~X⊃~(X•Y)∈TEO y ~Y⊃~(X•Y)∈TEO, por Ax+ y AxMP+ se infieren resultan +~X⊃+~(X•Y)∈TEO y +~Y⊃+~(X•Y)∈ TEO, es decir, ¬X⊃¬(X•Y)∈TEO y ¬Y⊃¬(X•Y)∈TEO, utilizando Ax1.5 y MP se concluye que (¬X∪¬Y)⊃ ¬(X•Y)∈TEO. □

**Proposición 3** (*Sustitución por equivalencia*). Sean X,Y,Z,W∈FOR.

a. Si X≡Y es un teorema entonces +X≡+Y, •X≡•Y, −X≡−Y, ¬X≡¬Y y *X≡*Y son teoremas.

b. Si F(Z)∈FOR es una fórmula en la cual figura Z y F(W) es el resultado de cambiar en F(Z) alguna ocurrencia de Z por W, entonces de Z≡W se infiere F(Z)≡F(W).

Prueba: Para la parte a, supóngase que, X≡Y, por EQ resulta (X⊃Y)•(Y⊃X), por la proposición 1 se infiere +((X⊃Y)•(Y⊃X)), por la proposición 2 se obtiene +(X⊃Y)•+(Y⊃X), utilizando I• y E• y AxMP+ se genera (+X⊃+Y)•(+Y⊃+X), finalmente, aplicando EQ se concluye +X≡+Y (Utilizando las definiciones de los operadores de verdad, también se siguen •X≡•Y, −X≡−Y, ¬X≡¬Y y *X≡*Y).

La parte b, se sigue de la parte a teniendo en cuenta SustEQ. □

## 4 SEMÁNTICA

Definición 4. (*Modelo, verdad*). ML es un *modelo de LD*, significa que, ML=(S, MA, , V), donde S es un conjunto no vacío de *mundos posibles,* MA es un mundo posible, llamado *mundo actual,* , es una relación binaria en MP, V es una *valuación* (función) de {MP, FOR} en {0, 1}.

La relación  satisface las siguientes restricciones:

RR. *Reflexividad.* (∀M∈S)MM.

RT. *Transitividad.* (∀M,N,F∈S)(MN y NF ⇒ MF).

RA. *Anti-simetría.* (∀M,N∈S)(MN y NM ⇒ M=N).

Sea **M**∈S y X,Y∈FOR y <u>a</u>∈FAA, la valuación **V** se define de la siguiente manera (V(M, X)=1 se abrevia como M(X)=1, y se dice que la fórmula X es *verdadera* en el mundo posible M):

1. V∼. M(∼X)=1 ⇔ M(X)=0     M(∼X)=0 ⇔ M(X)=1

2. V•. M(X•Y)=1 ⇔ M(X)=M(Y)=1     M(X•Y)=0 ⇔ M(X)=0 o M(Y)=0

3. V∪. M(X∪Y)=0 ⇔ M(X)=M(Y)=0     M(X∪Y)=1 ⇔ M(X)=1 o M(Y)=1

4. V⊃. M(X⊃Y)=0 ⇔ M(X)=1 y M(Y)=0     M(X⊃Y)=1 ⇔ si M(X)=1 entonces M(Y)=1

5. V≡. M(X≡Y)=1 ⇔ M(X)=M(Y)     M(X≡Y)=0 ⇔ M(X)≠M(Y)

6. V¬. M(¬X)=1 ⇔ (∀N∈MP)(M,N ⇒ N(X)=0)   M(¬X)=0 ⇔ (∃N∈MP)(M,N y N(X)=1)

7. V<u>a</u>.  M(<u>a</u>)=1 ⇔ (∀N∈MP)(M,N ⇒ N(<u>a</u>)=1)     M(<u>a</u>)=0 ⇔ (∃N∈MP)(M,N y N(<u>a</u>)=0)

Proposición 4 (*Caracterización semántica de los operadores de verdad*). En un modelo (S, M$_a$, , V) y Z∈FOR.

a. V+. V(M, +Z) = 1 ⇔ (∀N∈S)(MN ⇒ V(N, Z) = 1).

b. V⊗. V(M, ⊗Z) = 1 ⇔ (∃N∈S)(MN y V(N, Z) = 1).

c. V−. V(M, −Z) = 1 ⇔ (∃N∈S)(MN y V(N, Z) = 0).

d. V*. V(M, *Z) = 0 ⇔ (∃N∈S)(MN y V(N, Z) = 1) y (∃D∈S) (MD y V(D, Z) = 0).

Prueba. Supóngase que Z es una fórmula. De V¬ se tiene que, V(M, ¬∼Z) = 1 ⇔ (∀N∈S)(MN ⇒ V(N, ∼Z) = 0). Utilizando la definición de verdad alterna y la regla V∼ se infiere que V(M, +Z) = 1 ⇔ (∀N∈S)(MN ⇒ V(N, Z) = 1), por lo tanto, V+.

De V+ se tiene que, $V(M, +\sim Z) = 1 \Leftrightarrow (\forall N \in S)(MN \Rightarrow V(N, \sim Z) = 1)$. Por lo que, $V(M +\sim Z) = 0 \Leftrightarrow (\exists N \in S)(MN$ y $V(N, \sim Z) = 0)$, lo cual por V$\sim$ y la definición de $\neg$refutable significa que $V(M, \otimes Z) = 1 \Leftrightarrow (\exists N \in S)(MN$ y $V(N, Z) = 1)$, por lo tanto, V$\otimes$.

De V$\otimes$ se tiene que, $V(M, \otimes \sim Z) = 1 \Leftrightarrow (\exists N \in S)(MN$ y $V(N, \sim Z) = 1)$. Lo cual por V$\sim$ y las definiciones de $\neg$refutable y +refutable, implica que $V(M, -Z) = 1 \Leftrightarrow (\exists N \in S)(MN$ y $V(N, Z) = 0)$, por lo tanto, V$-$. Finalmente, observar que V* es consecuencia inmediata de V$\otimes$ y V$-$. □

## 5 VALIDEZ

Definición 5. (*Validez*). Cuando se tiene que $(\forall ML$ modelo$)(MA(X)=1)$, se dice que la fórmula X es **válida**, denotado X$\in$**VAL**. Dado un modelo $(S, M_a, , V)$, donde $M_a, M_2, \ldots, M_{t-2}, M_{t-1}, M_t$, son mundos posibles *diferentes* en S. Se dice que $C = M_a M_2 \ldots M_{t-2} M_{t-1} M_t$ es una *cadena de profundidad-t*, cuando se tienen $M_a M_2, \ldots, M_{t-2} M_{t-1}$ y $M_{t-1} M_t$.

Resulta entonces que X *no* es válida si y solamente si existe un modelo $M = (S, M_a, , V)$, en el cual $V(M_a, X) = 0$. Por lo que, si X no es válida, utilizando las reglas de las valuaciones, a partir de $V(M_a, X) = 0$, se construye un modelo $M = (S_a, M_a, , V)$ que refute la validez de la fórmula X, este modelo es llamado *modelo refutador*. Pero si la fórmula X es válida, entonces la construcción del *modelo refutador* fracasará, puesto que, en alguno de los mundos posibles (bien sea $M_a$ o un mundo generado por la aplicación de las reglas) del modelo en construcción se presentará una inconsistencia. Cuando fracasa la construcción del modelo refutador, entonces se genera una cadena de mundos posibles $C = M_a \ldots M_{k-1} M_k$ tal que *$M_k$ es inconsistente*, es decir, para alguna fórmula Z, $V(M_k, Z) = 1$ y $V(M_k, Z) = 0$. En este caso se dice que *la cadena C es inconsistente*.
En resumen, para probar la validez de una fórmula X, se supone que X no es válida, es decir, es falsa en el mundo actual $M_a$ de un modelo, y a partir de esta información se construye el modelo refutador. Si tal modelo no existe entonces se concluye que la fórmula X es válida.

Proposición 5 (*Preservación de la validez*). Para X$\in$FOR.
X$\in$VAL $\Rightarrow$ +X$\in$VAL.

Prueba. Supóngase que +X no es válida, por lo que existe un modelo refutador de +X, MO $= (S_{a+1}, M_{a+1}, , V)$ tal que $V(M_{a+1}, +X) = 0$, por la regla V+ resulta que existe un mundo posible $M_a$ tal que $M_{a+1} M_a$ y $V(M_a, X) = 0$. El (modelo MO se encuentra formado por cadenas consistentes de la forma $M_{a+1} M_a \ldots M_k$. A partir del modelo refutador MO de +X se construye un modelo refutador MO' de X de la siguiente manera: se elimina el mundo actual $M_{a+1}$ del conjunto $S_{a+1}$ obteniéndose el conjunto $S_a$, y se toma como mundo actual del modelo MO' el mundo $M_a$, en la relación de accesibilidad del modelo MO se eliminan las relaciones existentes entre el mundo actual $M_{a+1}$ y cualquier otro mundo obteniéndose la relación ', y del dominio de V se excluye el mundo actual $M_{a+1}$ obteniéndose V'. Como resultado se obtiene el modelo MO' $= (S_a, N_a, ', V')$, el cual por construcción se encuentra formado por cadenas consistentes $M_a \ldots S_k$, lo que significa que MO' es un modelo, y además en el mundo actual $M_a$ $V(M_a, X) = 0$, por lo tanto MO' es un modelo refutador de X, es decir X no es verdadera en el modelo MO', por lo que X no es válida. De lo anterior se concluye que si +X no es válida entonces X no es válida, es decir, X$\in$VAL $\Rightarrow$ +X$\in$VAL. □

Proposición 6 (*Validez de los axiomas*). Para X∈FOR.
X∈AxLD ⇒ X∈VAL.
Prueba. AxR. Si X es de la forma Z⊃⊗Z. Supóngase que X no fuese válida, entonces existiría un modelo tal que en el mundo actual $M_a$, V($M_a$, Z⊃⊗Z) = 0, lo cual según V⊃ significa V($M_a$, Z) = 1 y V($M_a$, ⊗Z) = 0, es decir V($M_a$, ~+~Z) = 0, resultando que V($M_a$, +~Z) = 1, utilizando la restricción RR se tiene $M_a M_a$, resultando que V($M_a$, ~Z) = 1, es decir V($M_a$, Z) = 0, lo cual no es el caso. Por lo tanto, Ax2.2∈VAL.

AxMP+. Si X es de la forma +(Y⊃Z)⊃(+Y⊃+Z). Si esta fórmula no fuese válida, entonces existiría un modelo tal que en el mundo actual $M_a$, V($M_a$, +(Y⊃Z)⊃(+Y⊃ +Z)) = 0, lo cual según la regla V⊃ significa V($M_a$, +(Y⊃Z)) = 1 y V($M_a$, +Y⊃+Z) = 0, y de nuevo por la misma regla se obtienen V($M_a$, +Y) = 1 y V($M_a$, +Z) = 0, de esta última por V+ se infiere la existencia de un mundo N, tal que $M_a$RN y V(N, Z) = 0, y como V($M_a$, +(Y⊃Z)) = 1 por V+ se infiere V(N, Y⊃Z) = 1, y como V($M_a$, +Y) = 1, por V+ se obtiene V(N, Y) = 1, y como ya se tiene V(N, Y⊃Z) = 1, por V⊃ se genera V(N, Z) = 1, pero esto es imposible. Por lo tanto, Ax2.1∈VAL.

AxCP. Si X es uno de los axiomas Ax1.1, ..., Ax1.13, utilizando las reglas V~, V•, V∪, V⊃ y V≡, y procediendo como es habitual para la validez del cálculo proposicional clásico en Caicedo (1990) y Hamilton (1981), se concluye que X∈VAL, es decir, Ax1.1, ..., Ax1.13∈VAL.

AxT. Sea X es de la forma ¬¬Z→⊗Z. Si X no fuese válida, entonces existiría un modelo tal que en el mundo actual $M_a$, V($M_a$, ¬¬Z⊃⊗Z) = 0, lo cual según V⊃ significa V($M_a$, ¬¬Z) = 1 y V($M_a$, ⊗Z) = 0, es decir V($M_a$, ~+~Z) = 0 y entonces V($M_a$, +~Z) = 1, utilizando la regla V− resulta que existe un mundo N, tal que $M_a$N, y en el cual V(N, ¬Z) = 0, por la regla V¬ resulta que existe un mundo K, tal que NK, y en el cual V(K, Z) = 1, pero por la restricción RT se obtiene $M_a$K, por lo que se infiere V(K, ~Z) = 1, es decir V(K, Z) = 0, lo cual es imposible, por lo tanto, ¬¬Z→⊗Z∈VAL, es decir, ¬Z→+Z∈VAL.

Sea X es de la forma $\underline{a}$→+$\underline{a}$, donde $\underline{a}$∈FAA. Si X no fuese válida, entonces existiría un modelo tal que en el mundo actual $M_a$, V($M_a$, $\underline{a}$→+$\underline{a}$) = 0, lo cual según V⊃ significa V($M_a$, $\underline{a}$) = 1 y V($M_a$, +$\underline{a}$) = 0; al tener V($M_a$, $\underline{a}$) = 1 por la regla V$\underline{a}$ resulta que para todo N, tal que $M_a$N, entonces V(N, $\underline{a}$) = 1; al tener V($M_a$, +$\underline{a}$) = 0 por V+ resulta que en todo mundo N, tal que $M_a$N, se infiere que V(N, $\underline{a}$) = 0, como $M_a M_a$, se obtienen V($M_a$, $\underline{a}$) = 1 y V($M_a$, $\underline{a}$) = 0, lo cual es imposible, por lo tanto, $\underline{a}$→+$\underline{a}$∈VAL. Al tener ¬Z→+Z∈VAL y $\underline{a}$→+$\underline{a}$∈VAL, se concluye que Ax2.3∈VAL.

Ax+. Si X es de la forma +Y, donde Y es uno de los anteriores axiomas, entonces resulta que, Y∈VAL, y por la proposición 5 se infiere que +Y∈VAL, por lo que Ax2.4∈VAL. □

Proposición 7 (*Validez*). Para X,Y∈FOR.
a. X∈VAL y X⊃Y∈VAL ⇒ Y∈VAL.
b. X∈TEO ⇒ X∈VAL.
Prueba. Supóngase que X∈VAL y X⊃Y∈VAL. Si Y no es válida, entonces existe un modelo tal que, en el mundo actual $M_a$, V($M_a$, Y) = 0. Como X∈VAL y X⊃Y∈VAL, entonces V($M_a$, X⊃Y) = 1 y V($M_a$, X) = 1, por la regla V⊃ de V($M_a$, Y) = 0 y V($M_a$, X⊃Y) = 1 resulta V($M_a$, X) = 0 lo cual es imposible. Por lo tanto Y∈VAL.
Para la parte b, supóngase que X∈TEO, se prueba que X∈VAL por inducción sobre la longitud L de la demostración de X.
**PB**. Paso Base L = 1. Si la longitud de la demostración de X es 1 entonces, X∈AxLD, lo cual por la proposición 6 significa que X∈VAL.

**PI**. Paso de inducción. Como hipótesis inductiva **HI**, se tiene que para cada fórmula Y, si Y∈TEO y la longitud de la demostración de Y es menor que L (donde L > 1) entonces Y∈VAL. Si X∈TEO y la longitud de la demostración de X es L entonces, X∈AxLD o X es consecuencia de aplicar MP en pasos anteriores de la demostración. En el primer caso se procede como en el caso base. En el segundo caso se tienen para alguna fórmula Y, demostraciones de Y y de Y⊃X, donde la longitud de ambas demostraciones es menor que L, utilizando la hipótesis inductiva se infiere que Y∈VAL y Y⊃X∈VAL, y por la parte a, resulta que X∈VAL. Utilizando el principio de inducción matemática, se ha probado la parte b. □

## 6 COMPLETITUD

Definición 6. (*Extensión consistente y completa*). Una *extensión* de LD, se obtiene alterando los axiomas de tal manera que, todos los teoremas de LD sigan siendo teoremas, y que el lenguaje de la extensión coincida con el lenguaje de LD. Una extensión es *consistente* si no existe ninguna X∈FOR tal que X como ∼X sean teoremas de la extensión. Un conjunto de fórmulas es *inconsistente* si de ellas se deriva una contradicción, es decir, si se deriva Z•∼Z para alguna Z∈FOR. Una extensión es *completa* si para toda X∈FOR, o bien X es teorema de la extensión o bien ∼X es teorema de la extensión. Para llegar a la prueba de completitud en la proposición 14, se sigue la metodología presentada en Henkin (1949).

Notación. **E∈EXT(LD)** significa que E es una extensión de LD. X∈**TEO-E** significa que X es un teorema de la extensión E.

Proposición 8 (*Extensión consistente*). Para X∈FOR.
a. LD es consistente.
b. Si **E**∈EXT(LD), X∉TEO-E y **E$_x$**∈EXT(LD) se obtiene añadiendo ∼X como nuevo axioma a E, entonces, E$_x$ es consistente.

Prueba. Parte a. Supóngase que LD no fuese consistente, por lo que debe existir X∈FOR tal que X,∼X∈TEO. Entonces por la proposición 7b, X,∼X∈VAL, pero esto es imposible, ya que si ∼X∈VAL, entonces para todo modelo (S, M$_a$, , V), se tienen V(M$_a$, ∼X) = 1, es decir, según V∼, V(M$_a$, X) = 0, por lo que X∉VAL, lo cual no es el caso. Por lo tanto, LD es consistente.

Parte b. Sea X∉TEO-E, y sea E$_x$ la extensión obtenida añadiendo ∼X como nuevo axioma a E. Supóngase que E$_x$ es inconsistente. Entonces, para alguna Z∈FOR, Z,∼Z∈TEO-E$_x$. Ahora bien, por Ax.9 se tiene que Z⊃(∼Z⊃X)∈TEO y por lo tanto de Z⊃(∼Z⊃X)∈TEO-E$_x$, aplicando dos veces MP se obtiene que X∈TEO-E$_x$. Pero E$_x$ tan sólo se diferencia de E en que tiene ∼X como axioma adicional, así que 'X es un teorema de E$_x$' es equivalente a 'X es un teorema de E a partir del conjunto {∼X}'. Por TD resulta que ∼X⊃X∈TEO-E, y por PR se infiere que X∈TEO-E, lo cual no es el caso. Por lo tanto, E$_x$ es consistente. □

Proposición 9 (*Extensión consistente y completa*).
Si E∈EXT(LD) es consistente de LD entonces existe E'∈EXT(E) que es consistente y completa.
Prueba. Sea X$_0$, X$_1$, X$_2$, . . . una enumeración de todas las fórmulas de LD. Se construye una sucesión J$_0$, J$_1$, J$_2$, . . . de extensiones de E como sigue: Sea J$_0$ = E. Si X$_0$∈TEO-J$_0$, sea J$_1$ = J$_0$. En caso contrario añádase ∼X$_0$ como nuevo axioma para obtener J$_1$ a partir de J$_0$. En general, dado t ≥ 1, para construir J$_t$ a partir de J$_{t-1}$, se procede de la siguiente manera: si X$_{t-1}$∈TEO- J$_{t-1}$, entonces J$_t$ = J$_{t-1}$, en caso contrario, sea J$_t$ la extensión de J$_{t-1}$ obtenida añadiendo ∼X$_{t-1}$ como nuevo axioma. La prueba se realiza por inducción matemática sobre t.
**PB**. Paso base. t=0. Como E es consistente por hipótesis, entonces J$_0$ es consistente.
**PI**. Paso inductivo. Hipótesis inductiva **HI**: J$_{t-1}$ es consistente con t ≥ 1. Por la proposición 8b, J$_t$ es consistente. Así pues, por el principio de inducción matemática, todo J$_t$ es consistente. Se define ahora J, como aquella extensión de E, la cual tiene como axiomas a aquellas fórmulas que son axiomas de al menos uno de los J$_t$.

Si J no es consistente, entonces existe X∈FOR tal que, X,~X∈TEO-J. Pero las demostraciones de X y ~X en J son sucesiones finitas de fórmulas, de modo que cada demostración solamente puede contener casos particulares de un número finito de axiomas de J. Por lo que, debe existir un t suficientemente grande, para que todos estos axiomas utilizados sean axiomas de $J_t$. Se deduce que X,~X∈TEO-J, lo cual es imposible ya que $J_t$ es consistente. Por lo tanto, J es consistente.

Para probar que J es completo, sea X∈FOR. X debe aparecer en la lista $X_0, X_1, X_2, \ldots$, supóngase que X es $X_k$. Si $X_k$∈TEO-$J_k$, entonces $X_k$∈TEO-J, puesto que J∈EXT($J_k$). Si $X_k$∉TEO-$J_k$, entonces de acuerdo con la construcción de $J_{k+1}$, ~$X_k$ es un axioma de $J_{k+1}$, con lo que ~$X_k$∈TEO-$J_{k+1}$, y entonces ~$X_k$∈TEO-J. Así, en todo caso se tiene que $X_k$∈TEO-J o ~$X_k$∈TEO-J, por lo que J es completo. □

Proposición 10 (*Consistencia subordinada*). Sean Y, $Z_1, \ldots, Z_k$∈FOR y $\underline{a}_{k+}, \ldots, \underline{a}_t$∈FAA.
Si {+$Z_1, \ldots,$ +$Z_k, \underline{a}_{k+1}, \ldots, \underline{a}_t$, ⊗Y} es consistente en LD entonces {$Z_1, \ldots, Z_k, \underline{a}_{k+1}, \ldots, \underline{a}_t$, Y} es consistente en LD.

Prueba. Si $\underline{a}$∈FAA, por AxT se tiene que $\underline{a}$⊃+$\underline{a}$, y AxR se tiene que +$\underline{a}$⊃$\underline{a}$, resultando que +$\underline{a}$≡$\underline{a}$.

Supóngase que {$Z_1, \ldots, Z_k, \underline{a}_{k+1}, \ldots, \underline{a}_t$, Y} es inconsistente en LD, por lo que existe una fórmula W∈FOR tal que, a partir de {$Z_1, \ldots, Z_k, \underline{a}_{k+1}, \ldots, \underline{a}_t$, Y} se infiere W•~W en LD, utilizando TD y Exp resulta que ($Z_1$• … • $Z_k$ • $\underline{a}_{k+1}$• … • $\underline{a}_t$ •Y)⊃(W•~W)∈TEO y por DI, en LD resulta ~($Z_1$• … • $Z_k$ • $\underline{a}_{k+1}$• … • $\underline{a}_t$ •Y)∈TEO, lo cual por N• e Imp significa, ($Z_1$ • … • $Z_k$ • $\underline{a}_{k+1}$• … • $\underline{a}_t$)⊃~Y∈TEO. Utilizando la proposición 1 resulta que +(($Z_1$•… • $Z_k$• $\underline{a}_{k+1}$• … • $\underline{a}_t$)⊃~Y)∈TEO, por AxMP+ se infiere +($Z_1$ • … • $Z_k$• $\underline{a}_{k+1}$• … • $\underline{a}_t$)⊃+~Y∈TEO, por la proposición 2a se obtiene (+$Z_1$ • … • +$Z_k$• +$\underline{a}_{k+1}$• … • +$\underline{a}_t$)⊃+~Y∈TEO, y como +$\underline{a}$≡$\underline{a}$, resulta que, (+$Z_1$ • … • +$Z_k$• $\underline{a}_{k+1}$• … • $\underline{a}_t$)⊃+~Y∈TEO, lo cual, por Imp, N• y la definición de ⊗, equivale a ~(+$Z_1$ • … • +$Z_k$ • $\underline{a}_{k+1}$• … • $\underline{a}_t$ ⊗Y)∈TEO, por lo que {+$Z_1, \ldots,$ +$Z_k, \underline{a}_{k+1}, \ldots, \underline{a}_t$, ⊗Y} es inconsistente en LD. Se ha probado que, {$Z_1, \ldots, Z_k, \underline{a}_{k+1}, \ldots, \underline{a}_t$, Y} inconsistente en LD implica que {+$Z_1, \ldots,$ +$Z_k, \underline{a}_{k+1}, \ldots, \underline{a}_t$, ⊗Y} inconsistente en LD, es decir, {+$Z_1, \ldots,$ +$Z_k, \underline{a}_{k+1}, \ldots, \underline{a}_t$, ⊗Y} consistente en LD implica {$Z_1, \ldots, Z_k, \underline{a}_{k+1}, \ldots, \underline{a}_t$, Y} consistente en LD. □

Definición 7. (*Subordinado*). Sean E,F∈EXT(LD) consistentes y completas. Se dice que F es **subordinado** de E si y solamente si existe Y∈FOR, tal que ⊗Y está en E, y además para cada Z∈FOR, tal que +Z está en E, y para cada $\underline{a}$∈FAA, tal que $\underline{a}$ está en E, se tiene que $\underline{a}$, Y y Z están en F.

Proposición 11 (*Extensión subordinada consistente y completa*). Para E∈EXT(LD), X∈FOR.
Si E es consistente y completa y ⊗X está en E, entonces existe F∈EXT(LD) consistente y completa tal que, X está en F y F subordinada de E.

Prueba. Sea X∈FOR, tal que ⊗X está en E. Sea $E_X$ = {X}∪{Z: +Z está en E}∪{$\underline{a}$: $\underline{a}$ está en E}, como E es consistente, entonces por la proposición 10, $E_X$ también es consistente. Al adicionar a $E_X$ los axiomas de LD y todas sus consecuencias, se obtiene una extensión de LD que incluye a $E_X$, utilizando la proposición 9 , se construye una extensión consistente y completa F de LD la cual incluye a $E_X$. Como X está en $E_X$, también está en F. Si +W está en E, por definición W está en $E_X$, por lo que W está en F. Si $\underline{a}$ está en E, por definición $\underline{a}$ está en $E_X$, por lo que W está en F. Por lo que, F es subordinado de E. □

Proposición 12 (*Propiedades de la subordinación*). Para E,F,G∈EXT(LD) consistentes y completas.
a. (*Transitividad*). (F es subordinado de E y G es subordinado de F) ⇒ G es subordinado de E.
b. (*Reflexividad*). F es subordinado de F.
c. (*Anti-simetría*). (F es subordinado de G y G es subordinado de F) ⇒ F=G.
Prueba. Parte a, supóngase que G es subordinado de F y F es subordinado de E. Como G es subordinado de F entonces existe en F, ⊗Z∈FOR tal que Z está en G. Si ⊗Z no está en E, entonces al ser una extensión completa,

∼⊗Z si debe estarlo, además por AxT se tiene que ⎯⎯¬Z⊃⊗Z está en E, por lo que ∼⎯⎯¬Z, es decir ++∼Z también está en E, y al ser F subordinado de E resulta que +∼Z está en F, lo cual significa que ∼⊗Z está en F, pero esto es imposible ya que F es consistente. Por lo que, ⊗Z está en E.

Sea W∈FOR, tal que +W está en E, es decir ∼⊗∼W está en E, utilizando AxT ⎯⎯¬∼W⊃⊗∼W, resulta que ∼⎯⎯¬∼W está en E, por lo que ++W está en E, y como F es subordinado de E, se infiere que +W está en F, como además, G es subordinado de F, entonces W está en G.

Sea $\underline{a}$∈FAA, tal que $\underline{a}$ está en E, utilizando AxT $\underline{a}$⊃+$\underline{a}$, resulta que +$\underline{a}$W está en E, por lo que ++$\underline{a}$ está en E, y como F es subordinado de E, se infiere que +$\underline{a}$ está en F, como G es subordinado de F, entonces $\underline{a}$ está en G.

En resumen, existe ⊗Z en E tal que para cada fórmula +W en E se tiene que Z y W están en G, y por lo tanto, G es subordinado de E.

Parte b. Sea X el axioma Ax1.1, por lo que en LD se tiene X, y como por AxR se tiene X⊃⊗X, resulta ⊗X, por lo que X y ⊗X están en F. Supóngase que +W está en F, por AxR en F se tiene ∼W⊃⊗∼W, es decir ∼⊗∼W⊃W, o sea +W→W, resultando que W también está en F. Por lo tanto, F subordinada de F.

Parte c. Sean F es subordinado de G y G es subordinado de F. Sea X∈FOR. Si X∈F entonces +X∈F y como G es subordinado de F entonces X∈G. Si X∈G entonces +X∈G y como F es subordinado de G entonces X∈F. Se tiene entonces que (∀X∈FOR)(X∈F ⇔ X∈G), lo cual significa que F=G. □

Proposición 13 (*Construcción de un modelo*).

Si E'∈EXT(LD) es consistente, entonces existe un modelo en el cual todo X∈TEO-E' es verdadero.

Prueba. Se define el modelo (S, ME$_n$, , V) de la siguiente manera: sean E, F, G, ..., extensiones consistentes y completas de E' (E$_a$ la inicial y las demás subordinadas), presentadas en las proposiciones 9 y 11. A cada extensión F, se le asocia un mundo posible MF, sean S el conjunto de tales mundos posibles y ME$_a$ el mundo actual.

La relación de accesibilidad se construye así: MF$_t$MG$_t$ si y solamente si G$_t$ es subordinado de F$_t$.

Para cada MF$_k$ en S y para cada X∈FOR, V(MF$_k$, X) = 1 si X está en F$_k$ y V(MF$_k$, X) = 0 si ∼X está en F$_k$, donde F$_k$ es la extensión consistente y completa asociada a MF$_k$. Nótese que V es funcional, por ser F$_k$ consistente y completa. Ahora bien, ya que F$_k$ es consistente, entonces V(MF$_k$, X) ≠ V(MF$_k$, ∼X) y por lo tanto, V(MF$_k$, X) = 1 ⇔ V(MF$_k$, ∼X) = 0, por lo que se satisface la definición V∼. Para afirmar que M es un modelo, se debe garantizar que, para cada uno de los conectivos, V satisface la definición de valuación.

Para el caso del condicional X⊃Y. Utilizando N⊃, I• y E•, se tiene la siguiente cadena de equivalencias: V(MF$_k$, X⊃Y) = 0, es decir ∼(X⊃Y) está en F$_k$, o sea que X•∼Y está en F$_k$, resultando que X y ∼Y están en F$_k$, lo cual significa que V(MF$_k$, X) = 1 y V(MF$_k$, Y) = 0, por lo que se satisface la definición V⊃.

Para el caso de la conjunción X•Y. Utilizando I• y E•, se tiene la siguiente cadena de equivalencias: V(MF$_k$, X•Y) = 1, es decir X•Y está en F$_k$, por lo que X y Y están en F$_k$, lo cual significa que V(MF$_k$, X) = 1 y V(MF$_k$, Y) = 1, por lo que se satisface la definición V•.

Para el caso de la disyunción X∪Y. Utilizando N∪, I• y E•, se tiene la siguiente cadena de equivalencias: V(MF$_k$, X∪Y) = 0, es decir ∼(X∪Y) está en F$_k$, o sea que ∼X•∼Y está en F$_k$, de donde ∼X y ∼Y están en F$_k$, es decir V(MF$_k$, X) = 0 y V(MF$_k$, Y) = 0, por lo que se satisface la definición V∪.

Para el caso del bicondicional X≡Y. Utilizando EQ, V∪, V• y V∼ (ya probadas) e I• y E•, se tiene la siguiente secuencia de equivalencias: V(MF$_k$, X≡Y) = 1, es decir X≡Y está en F$_k$, por lo que (X•Y)∪(∼X•∼Y) está en F$_k$, lo que significa V(MF$_k$, (X•Y)∪(∼X•∼Y)) = 1, o de otra forma V(MF$_k$, X•Y) = 1 o V(MF$_k$, ∼X•∼Y) = 1, es decir, V(MF$_k$, X) = V(MF$_k$, Y) = 1 o V(MF$_k$, X) = V(MF$_k$, Y) = 0, o dicho de otra manera V(MF$_k$, X) = V(MF$_k$, Y), por lo que se satisface la definición V≡.

Para el caso de la regla V¬. MF$_p$ es un mundo asociado a F$_p$, MG$_p$ es un mundo asociado a G$_p$, y Z∈FOR. Supóngase que V(MF$_p$, ¬Z) = 1, por lo que ¬Z, es decir +∼Z está en F$_p$. Si MF$_p$RMG$_p$, entonces G$_p$ es subordinada de $F_p$ y entonces $\sim Z$ está en $G_p$, resultando que $V(MG_p, Z) = 0$. Se ha probado de esta manera que $V(MF_p, \neg Z) = 1 \Rightarrow (\forall MG_p \in S)(MF_p RMG_p \Rightarrow V(MG_p, Z) = 0)$.

Para probar la recíproca, supóngase que $(\forall MG_p \in S)(MF_p RMG_p \Rightarrow V(MG_p, Z) = 0)$. Si $V(MF_p, \neg Z) = 0$, entonces al ser $MF_p$ el mundo asociado a la extensión consistente y completa $F_p$ resulta que $\sim\neg Z$ está en $F_p$, por lo que $\sim+\sim Z$, es decir $\otimes Z$ está en $F_p$. Por la proposición 11 existe una extensión consistente y completa $G_p$ subordinada de $F_p$ tal que $Z$ está en $G_p$. Como $MG_p$ es el mundo asociado a $G_p$, entonces $MF_p RMG_p$, lo cual, por el supuesto inicial implica $V(MG_p, Z) = 0$, es decir $\sim Z$ está en $G_p$, resultando que $G_p$ es inconsistente, lo cual no es el caso. Por lo tanto, $V(MF_p, \neg Z) = 1$. Se ha probado de esta manera que $(\forall MG_p \in S)(MF_p RMG_p \Rightarrow V(MG_p, Z) = 0) \Rightarrow V(MF_p, \neg Z) = 1$. Resultando que se satisface V¬.

Para el caso de la regla V$\underline{a}$. $MF_p$ es un mundo asociado a $F_p$, $MG_p$ es un mundo asociado a $G_p$, y $\underline{a} \in$ **FAA**. Supóngase que $V(MF_p, \underline{a}) = 1$, por lo que $\underline{a}$ está en $F_p$ y como AxT $\underline{a} \supset +\underline{a}$ está en $F_p$, entonces por MP también $+\underline{a}$ está en $F_p$. Si $MF_p RMG_p$, entonces $G_p$ es subordinada de $F_p$ y entonces $\underline{a}$ está en $G_p$, resultando que $V(MG_p, \underline{a}) = 1$. Se ha probado de esta manera que $V(MF_p, \underline{a}) = 1 \Rightarrow (\forall MG_p \in S)(MF_p RMG_p \Rightarrow V(MG_p, \underline{a}) = 1)$.

Para probar la recíproca, supóngase que $(\forall MG_p \in S)(MF_p RMG_p \Rightarrow V(MG_p, \underline{a}) = 1)$ donde $\underline{a} \in$ **FAA**. Si $V(MF_p, \underline{a}) = 0$, entonces al ser $MF_p$ el mundo asociado a la extensión consistente y completa $F_p$ resulta que $\sim \underline{a}$ está en $F_p$, utilizando AxR $+\underline{a} \supset \underline{a}$, se sigue $\sim +\underline{a}$, es decir $\otimes \sim \underline{a}$ está en $F_p$. Por la proposición 11 existe una extensión consistente y completa $G_p$ subordinada de $F_p$ tal que $\sim \underline{a}$ está en $G_p$. Como $MG_p$ es el mundo asociado a $G_p$, entonces $MF_p RMG_p$, lo cual, por el supuesto inicial implica $V(MG_p, \underline{a}) = 1$, es decir $\underline{a}$ está en $G_p$, resultando que $G_p$ es inconsistente, lo cual no es el caso. Por lo tanto, $V(MF_p, \underline{a}) = 1$. Se ha probado de esta manera que $(\forall MG_p \in S)(MF_p RMG_p \Rightarrow V(MG_p, \underline{a}) = 0) \Rightarrow V(MF_p, \underline{a}) = 1$. Como ya se probó la recíproca, resulta que se satisface V$\underline{a}$.

Con base en el análisis anterior, y teniendo en cuenta que las reglas RT, RR y RA, se encuentran garantizadas por las proposiciones 11 y 12 y la forma en que se construye el modelo, se concluye finalmente que V es una valuación, y por lo tanto, M es un modelo.

Para finalizar la prueba, sea X un teorema de E', por lo que X está en E'. Por lo tanto, utilizando la definición de V resulta que $V(ME_a, X) = 1$, es decir, X es verdadera en el modelo $M = (S, ME_a, , V)$. □

Proposición 14 (*Completitud de LD*). Para $X \in$ FOR.
$X \in$ VAL $\Rightarrow X \in$ TEO.

Prueba. Si $X \notin$ TEO, entonces, por la proposición 8b, la extensión E', obtenida añadiendo $\sim X$ como nuevo axioma, es consistente. Así pues, según la proposición 13, existe un modelo M tal que todo teorema de E' es verdadero en M, y como $\sim X \in$ TEO-E', entonces $\sim X$ es verdadero en M, es decir, X es falso en M, y por lo tanto, $X \notin$ VAL. Se ha probado que, $X \notin$ TEO $\Rightarrow X \notin$ VAL, es decir, $X \in$ VAL $\Rightarrow X \in$ TEO. □

Proposición 15 (*Caracterización semántico deductiva de LD*). Para $X \in$ FOR.
$X \in$ VAL $\Leftrightarrow X \in$ TEO.
Prueba: Consecuencia de las proposiciones 7b y 14. □

# 7 GRÁFICOS EXISTENCIALES GAMMA-LD

En esta sección se resumen los gráficos existenciales Gamma-LD prenetados en Sierra-Aristizabal (2021). Se toma como punto de partida la *Hoja de aserción* H donde se dibujan los gráficos existenciales. El conjunto Alfa-LC de *gráficos existenciales clásicos*, se define de la siguiente manera: 1) a∈FAC ⇒ a∈Alfa-LC. 2) λ∈Alfa-LC (λ es el gráfico vacío). 3) X∈Alfa-LC ⇒ (X)∈Alfa-LC. 4) X,Y∈Alfa-LC ⇒ XY∈Alfa-LC. 5) Sólo 1 a 4 determinan Alfa-LC.

El conjunto Alfa-LI de *gráficos existenciales intuicionistas*, se define de la siguiente manera: 1) a∈FAA ⇒ a∈Alfa-LI. 2) λ∈Alfa-LI (λ es el gráfico vacío). 3) X∈Alfa-LI ⇒ [X]∈Alfa-LI. 4) X,Y∈Alfa-LI ⇒ XY, [X (Y)], [(X) (Y)]∈Alfa-LI. 5) Sólo 1 a 4 determinan Alfa-LI.

El conjunto Gamma-LD de *gráficos existenciales LD* se define de la siguiente manera: 1) X∈Alfa-LC ⇒ X∈Gamma-LD. 2) λ∈Gamma-LD (λ es el gráfico vacío). 3) X∈Alfa-LI ⇒ X∈Gamma-LD. 4) X∈Gamma-LD ⇒ (X),[X]∈Gamma-LD. 5) X,Y∈Gamma-LD ⇒ XY,[X (Y)],[(X) (Y)]∈Gamma-LD. 6) Sólo 1 a 5 determinan Gamma-LD.

El conjunto GA de gráficos existenciales alternos se define de la siguiente manera: X∈GA ⇔ X=[Y] con Y∈Gamma-LD o X=a donde a∈FAA. Notación. Si X∈GAE entonces $\underline{X}$ indica que X∈GA.

**Reglas de transformación y Validez.** Cuando se tiene el gráfico existencial (X), se dice X está rodeado por un *corte clásico* y cuando se tiene el gráfico existencial [X], se dice X está rodeado por un *corte paracompleto*. Se dice que un gráfico existencial X se encuentra en una *región par*, denotado X∈RP, si X se encuentra rodeado por un número par de cortes (clásicos y/o paracompletos). X se encuentra en una *región impar*, denotado X∈RI, si X se encuentra rodeado por un número impar de cortes (clásicos y/o paracompletos). X se encuentra en una *región clásica*, denotado X∈RC, si X no se encuentra rodeado por cortes paracompletos. X se encuentra en una *región alterna*, denotado X∈RA, si X se encuentra rodeado por al menos un corte paracompleto.

Notación. $X \overset{Re}{\Rightarrow} Z$ significa que la *regla de transformación Re*, aplicada a un gráfico X permite inferir un nuevo gráfico Z. Si también Z se transforma en X mediante *Re* entonces se escribe $X \overset{Re}{\Leftrightarrow} Z$.

X>>Z significa que X se transforma en Z utilizando un número finito de reglas de transformación. Para X∈Gamma-LD, se dice que X es *válido*, denotado X∈GEV, si λ>>X.

Definición . (*Definición de λ*). Def-λ. Definición de λ. $\lambda \overset{Def-\lambda}{\Longleftrightarrow}$ ' ' (gráfico vacío).

**Reglas de transformación primitivas en Gamma-LD.** Notación. Un subíndice ₁ indica que la región donde se encuentra el gráfico es *impar*, un subíndice ₂ indica que la región es par.

1) Rλ: Regla lambda, la hoja de aserción es un gráfico existencial válido, λ∈GEV. 2) B: Borrado, $XY_2 \overset{B}{\Rightarrow} X_2$, $XY_2 \overset{B}{\Rightarrow} Y_2$. 3) E: Escritura, $X_1 \overset{E}{\Rightarrow} XY_1$, $X_1 \overset{E}{\Rightarrow} YX_1$. 4) DCC: Doble corte clásico, $X \overset{DCC}{\Leftrightarrow} ((X))$. 5) CC: Cambio de corte, $[X]_2 \overset{CC}{\Rightarrow} (X)_2$, $(X)_1 \overset{CC}{\Rightarrow} [X]_1$. 6) DCMGEV: Doble corte mixto para Gráficos existenciales válidos, X∈GEV ⇒ $X \overset{DCMGEV}{\Longleftrightarrow} [(X)]$. 7) DCMF: Doble corte mixto fuerte, $\underline{X}_2 \overset{DCMF}{\Longrightarrow} [(\underline{X}_2)]$, $[(\underline{X}_1)] \overset{DCMF}{\Longrightarrow} \underline{X}_1$. 8) I: Iteración, $X \overset{I}{\Rightarrow} XX$. 9) D: Desiteración, $XX \overset{D}{\Rightarrow} X$. 10) IC: Iteración clásica, $X(…(Y)…) \overset{IC}{\Rightarrow} X(…(XY)…)$. 11) DC: Desiteración clásica, $X(…(XY)…) \overset{DC}{\Rightarrow} X(…(Y)…)$. 12) IF: Iteración fuerte, $\underline{X}/…/Y\backslash…\backslash \overset{IF}{\Rightarrow} \underline{X}/…/\underline{X}Y\backslash…\backslash$ Donde /Z\ significa (Z) o [Z]. 13) DF: Desiteración fuerte, $\underline{X}/…/\underline{X}Y\backslash…\backslash \overset{DF}{\Rightarrow} \underline{X}/…/Y\backslash…\backslash$, 14) Sólo aplican las reglas 1 a 13 y sus reglas derivadas. Notación. Las reglas de la forma $X \overset{R}{\Leftrightarrow} Y$ indican que $X \overset{R}{\Rightarrow} Y$ y $Y \overset{R}{\Rightarrow} X$ valen tanto en regiones pares como en regiones impares.

Proposición. (*Reglas de transformación derivadas*). 1) DCCλ: Doble corte clásico, λ,((λ)),(( ))∈GEV. 2) DCM: Doble corte mixto, $[(X_2)] \overset{DCM}{\Longrightarrow} X_2$, $X_1 \overset{DCM}{\Longrightarrow} [(X)]_1$. 3) DCMλ: Doble corte mixto, λ,[(λ)],[( )]∈GEV. 4) CCE: Cambio de corte especifico, $[X (Y)]_2 \overset{CCE}{\Longrightarrow} [X [Y]]_2$, $[X [Y]]_1 \overset{CCE}{\Longrightarrow} [X (Y)]_1$. 5) DCMF.1: Doble corte mixto fuerte, $[(\underline{X})] \overset{DCMF.1}{\Longleftrightarrow} \underline{X}$. 6) TCM: Triple corte mixto, $[([X])] \overset{TCM}{\Longleftrightarrow} [X]$. 7) DCAF: Doble corte alterno fuerte, $\underline{X}_2 \overset{DCAF}{\Longrightarrow} [[\underline{X}_2]]$, $[[\underline{X}_1]] \overset{DCAF}{\Longrightarrow} \underline{X}_1$. 8) TCA: Triple corte alterno $[X]_2 \overset{TCA}{\Longrightarrow} [[[X]_2]]$, $[[[X]_1]] \overset{TCA}{\Longrightarrow} [X]_1$. 9)

TCAF: Triple corte alterno, fuerte, $[[\underline{[X]}_2]] \overset{TCAF}{\Longrightarrow} [\underline{X}]_2$, $[\underline{X}]_1 \overset{TCAF}{\Longrightarrow} [[\underline{[X]}_1]]$. 10) TCAF.1: Triple corte alterno fuerte, $[[\underline{[X]}]] \overset{TCAF.1}{\Longleftrightarrow} [\underline{X}]$. 11) CCA: Cuádruple corte alterno, $[[[[X]]]] \overset{CCA}{\Longleftrightarrow} [[X]]$.

Proposición. (*Reversión de las reglas de transformación*). Para X,Y∈Gamma-LD. a) ($\forall$R∈RTRA) $[X_2 \overset{R}{\Longrightarrow} Y](\exists R'\in RTRA)[Y_1 \overset{R'}{\Longrightarrow} X]$. b) ($\forall$R∈RTRA)$\{[\ X_2 \overset{R}{\Longrightarrow} Y \Rightarrow Y_2]$ y $[X_1 \overset{R}{\Longrightarrow} Y \Rightarrow Y_1]\}$

Proposición. (*Principio de contraposición*). Para $X_0, X_n\in$Gamma-LD. $X_0 >> X_n \Rightarrow \{(X_0 \in RP \Rightarrow X_0 >> X_n)$ y $(X_n \in RI \Rightarrow X_n >> X_0)\}$

Proposición. (*Teoremas de deducción gráfico y fuerte en Gamma-LD*). Para X,Y∈Gamma-LD. a) TDG. $X>>Y \Rightarrow (X(Y))$. b) TDGF. $\underline{X}>>\underline{Y} \Rightarrow [\underline{X}(\underline{Y})]$.

Proposición. (*Teorema de demostración indirecta en Gamma-LD*). Para X∈Gamma-LD. a) $X>>(\ ) \Rightarrow (X)$. b) $(X)>>(\ ) \Rightarrow X$. c) $\underline{X}>>[\ ] \Rightarrow [\underline{X}]$.

**Traducción de FOR a Gamma-LD.** Definición. (*Traducción de fórmulas a gráficos existenciales*). La función ( )' de FOR en Gamma-LD, se define de la siguiente manera (donde a∈**FAT** y X,Y): 1) a' = a. 2) {~X}' = (X'). 3) {X•Y}' = X'Y'. 4) {X⊃Y}' = (X'(Y')). 5) {X≡Y}' = {X⊃Y}' {Y⊃X}' = (X'(Y')) (Y'(X')). 6) {X∪Y}' = ((X')(Y')). 7) Ax1.1' = λ. 8) {¬X}' = [X']. 9) {+X}' = [(X')]. 10) {X→Y}' = [X'(Y')]. 11) {X∨Y}' = [(X')(Y')]. 12) {X∧Y}' = [(X' Y')]. 13) {X↔Y}' = [X'(Y')] [Y'(X')].

El resultado principal probado en Sierra-Aristizabal (2021) es el siguiente:

Proposición 16. (*Los teoremas de LD son exactamente los gráficos válidos en Gamma-LD*). ($\forall Z \in$FOR) (Z∈TEO $\Leftrightarrow$ Z'∈GEV)

Proposición 17 (*Caracterización semántico de los gráficos existenciales Gamma-LD*). Para X∈FOR. X∈VAL $\Leftrightarrow$ X∈TEO.

Prueba: Consecuencia directa de las proposiciones 15 y 16. □

# 8 CONCLUSIONES

Conclusión 1. (*Restricción de FOR a FC*). Si los modelos de LD se restringen al conjunto FC de fórmulas clásicas, resulta que el conjunto de fórmulas válidas, en el lenguaje {~, ⊃, •, ∪, ≡}, coincide con el conjunto de teoremas del cálculo proposicional clásico Van Dale (2013). Por lo que, los teoremas del cálculo proposicional clásico son válidos en LD. □

Conclusión 2. (*Restricción de FOR a FI*). Si los modelos de LD se restringen al conjunto FI de fórmulas intuicionistas, en el lenguaje {¬, ∧, ∨ →, ↔}, resulta que el conjunto de fórmulas válidas coincide con el conjunto de teoremas del cálculo proposicional intuicionista Van Dale (2013). Por lo que, los teoremas del cálculo proposicional intuicionista son válidos en LD.
Prueba: Basta con verificar el significado de los conectivos intuicionistas en LD.

V¬. V(M, ¬Z) = 1 $\Leftrightarrow$ ($\forall$N∈S)(MN $\Rightarrow$ V(N, X⊃Z) = 0).

V∧. V(M, X∧Z) = 1 $\Leftrightarrow$ V(M, +(X•Z) = 1 $\Leftrightarrow$ ($\forall$N∈S)(MN $\Rightarrow$ V(N, X•Z) = 1).

V∨. V(M, X∨Z) = 1 $\Leftrightarrow$ V(M, +(X∪Z) = 1 $\Leftrightarrow$ ($\forall$N∈S)(MN $\Rightarrow$ V(N, X∪Z) = 1).

V→. V(M, X→Z) = 1 $\Leftrightarrow$ V(M, +(X⊃Z) = 1 $\Leftrightarrow$ ($\forall$N∈S)(MN $\Rightarrow$ V(N, X⊃Z) = 1).

V↔. V(M, X↔Z) = 1 $\Leftrightarrow$ V(M, +(X≡Z) = 1 $\Leftrightarrow$ ($\forall$N∈S)(MN $\Rightarrow$ V(N, X≡Z) = 1). □

Conclusión 3. (LD *generaliza a LI, a LC y a S4*). Como consecuencia de la dos conclusiones anteriores se tiene que la Lógica Doble es una generalización de la Lógica Proposicional Clásica LC, de la Lógica Proposicional Intuicionista LI y de la Lógica Proposicional Modal S4. ☐

Platón en uno de sus diálogos, *Crátilo*, define la verdad como "El discurso, que dice las cosas como son, es verdadero; y el que las dice como no son, es falso" Platón (1983). En el libro IV de la Metafísica, Aristóteles define el concepto de verdad de la siguiente manera "decir de *lo que es* que es, y de *lo que no es* que no es, es *lo verdadero*; decir de *lo que es* que no es, y de *lo que no es* que es, es *lo falso*" Aristóteles (1998). Con el sistema LD se puede modelar esta definición interpretando +X como 'X es verdadero', ¬X como 'X es falso', *X como 'decir X', X como 'X es' y ~X como 'X no es'. La definición de *verdad aristotélica* estaría codificada por la fórmula (A•*A)≡+A, y la definición de *falsedad aristotélica* estaría codificada por la fórmula (~A•*A)≡¬A.

Conclusión 4. (*Verdad y falsedad aristotélica*). Para X∈FOR.
a. (X•*X)≡+X∈TEO.
b. (~X•*X)≡¬X∈TEO.
Prueba: Parte a. Supóngase X•*X. Por E• se obtienen X y *X, por AxR ¬X⊃~X, es decir X⊃~¬X, se infiere ~¬X, al tener *X, resulta +X∪¬X. Aplicando SD en estos dos resultados se obtiene +X. Por TD se ha probado (X•*X)⊃+X.
Por Ax0.3 se tiene +X⊃(+X∪¬X), lo cual significa +X⊃*X, de AxR se infiere +X→X, de estos dos resultados, utilizando Ax0.8 se obtiene +X⊃(X•*X).
Como se tienen (X•*X)⊃+X y +X⊃(X•*X), utilizando Ax0.11 se infiere (X•*X) ≡ +X.
Para la parte b. Supóngase ~X•*X. Por E• se obtienen ~X y *X, utilizando AxR +X⊃X se infiere ~+X, al tener *X, resulta +X∪¬X, aplicando SD en estos dos resultados se obtiene ¬X. Por TD se ha probado (~X•*X)⊃¬X.
Por Ax0.4 se tiene ¬X⊃(+X∪¬X), lo cual significa ¬X⊃*X y por AxR se infiere ¬X⊃~X, de estos resultados, utilizando Ax0.8 se obtiene ¬X⊃(~X•*X).
Como se tienen (~X•*X)⊃¬X y ¬X⊃(~X•*X), utilizando Ax0.11 se infiere (~X•*X) ≡ ¬X. ☐

Muchas *paradojas lógicas* involucran los conceptos de verdad o falsedad, por ejemplo, la siguiente variante de la paradoja del mentiroso Bochenski (1976) y Smullyan (1997). Considérese la situación en la cual se tiene una oración que dice:

> Esta oración es falsa

Cuando se identifican en la lógica clásica, ser el caso con verdadero (Z) y no ser el caso con falso (~Z), entonces se tiene la paradoja: si la oración es el caso (Z) entonces resulta que también es falsa (~Z), y si la oración es falsa (~Z) entonces resulta que es el caso (Z). Es decir, (Z⊃~Z)•(~Z⊃Z). Resultando Z•~Z. Como la lógica clásica no soporta las contradicciones, resulta inútil en este caso.

Cuando se identifican en la lógica intuicionista, ser el caso como verdadero (Z), y no ser el caso con falso (¬Z), entonces se tiene la paradoja: si la oración es el caso (Z) entonces resulta que también es falsa (¬Z), y si la oración es falsa (¬Z) entonces resulta que es el caso (Z). Es decir, (Z→¬Z)•(¬Z→Z). Resultando Z∧¬Z. Como la lógica intuicionista no soporta las contradicciones, resulta inútil en este caso.

Cuando se identifican en la lógica intuicionista, ser el caso (Z), la doble negación como verdadero (¬¬Z), y no ser el caso con falso (¬Z), entonces se tiene la paradoja: si la oración es verdadera (¬¬Z) entonces resulta que también es falsa (¬Z), y si la oración es falsa (¬Z) entonces resulta que es verdadera (¬¬Z). Es decir, (¬¬Z→¬Z)•(¬Z→¬¬Z). Resultando ¬Z∧¬¬Z. Como la lógica intuicionista no soporta las contradicciones, resulta inútil en este caso.

Comentario 1. (*Las fórmulas inconsistentes no tienen modelos*) Para X∈FOR,

Si existe un modelo (S, M, R, V) de LD tal que V(M, X)=1 entonces X es consistente con LD.

Prueba. Parte a. Sea X inconsistente con LD, por lo que X≡(Y•~Y) para algún Y∈FOR. Si hay un modelo tal que V(M, X)=1, entonces por V≡ se sigue que V(M, Y•~Y )=1, lo cual por V• y V~ implica que V(M, Y)=1 y V(M, Y)=0, es decir 1=0, lo cual no es el caso, y por lo tanto, no existe tal modelo. Se ha probado que, Si X es inconsistente con LD entonces no existe un modelo (S, M, , V) de LD tal que V(M, X)=1. Lo cual equivale a:

si existe un modelo (S, M, , V) de LD tal que V(M, X)=1 entonces X es consistente con LD. □

Conclusión 5. (*Solución a una paradoja*). En la Lógica Doble LD, representando *lo que dice la oración es el caso* como Z, *lo que dice la oración no es el caso* como ~Z, *la oración Z es falsa* como ¬Z, y *la oración Z es verdadera* como +Z, entonces la situación (se tiene una oración que dice: *esta oración es falsa*) queda representada por la fórmula Z≡¬Z, la cual en LD no genera contradicciones y se pueden obtener conclusiones válidas.

Prueba: Si la oración Z es verdadera (+Z), por AxR resulta que lo que dice es el caso (Z), y como la oración dice que Z es falsa (Z≡¬Z), se sigue que Z es falsa (¬Z), y utilizando AxR se infiere que Z no es el caso (~Z), y de nuevo por AxR, Z no puede ser verdadera (~+Z), se obtiene de esta manera una contradicción (+Z•~+Z), se concluye por demostración indirecta, que **Z no es verdadera** (~+Z).

Ahora, si Z es falsa (¬Z) entonces por AxR resulta lo que dice Z no es el caso (~Z), y como la oración dice que Z es falsa (Z≡¬Z), se sigue que Z no es falsa (~¬Z), se obtiene una contradicción (¬Z•~¬Z), por demostración indirecta, se concluye que **Z no es falsa** (~¬Z).

Si Z es el caso (Z), como la oración dice que Z es falsa (Z≡¬Z), se sigue que Z es falsa (¬Z), pero ya se probó que Z no es falsa (~¬Z). Por lo tanto, **Z no es el caso** (~Z).

Se ha probado Z ni es verdadera ni es falsa (~+Z•~¬Z), lo cual significa que **Z no está bien fundada** (~*Z), este hecho implica que **la situación no está bien fundada** (~*(Z≡¬Z)), la prueba se deja al lector.

De acuerdo con el comentario 1, para garantizar que no hay paradoja, es decir que Z≡¬Z no genera contradicción, basta verificar que Z≡¬Z tiene un modelo. Considere el modelo ({MA, M1}, MA, , V) tal que, MAMA, M1M1, MAM1, V(MA, ¬Z) = V(MA, +Z) = V(MA, Z) = 0 y V(M1,Z)=1. Por lo tanto, V(MA, Z≡¬Z) =1. Es decir, Z≡¬Z tiene un modelo, y por lo tanto es consistente, y no hay paradoja. □

Comentario 2. (*Dos fórmulas invalidas en LD*) Para X∈FOR.
a. ~+X⊃¬+X∉TEO.
b. ¬+X⊃ ¬X∉TEO.

Prueba. Parte a. Considérese el modelo ML =(S, MA, , V), tal que S={MA, M1, M2}, MAM1, MAM2, MAMA, M1M1, M2M2, M1(X)=0, M2(X)=1. Como M1(X)=0 y MAM1 entonces por V+ resulta que MA(+X)=0, por V~ se infiere que MA(~+X)=1. Como M2(X)=1 y sólo M2M2, entonces por V+ se obtiene M2(+X)=1, por V¬ se deduce MA(¬+X)=0. En consecuencia, por V⊃ se puede asegurar que MA(~+X⊃¬+X)=0, lo cual significa que ~+X⊃¬+X no es válida, y según la proposición 7b, ~+X⊃¬+X∉TEO.

Parte b. Considérese el modelo ML =(S, MA, , V), tal que S={MA, M1, M2}, MAM1, MAM2, MAMA, M1M1, M2M2, M2M1, M1(X)=0, M2(X)=1.

Como M1(X)=0, MAM1, M1M1 y M2M1 por V+ se sigue M1(+X)=M2(+X)=MA(+X)=0, resultando por V¬ que MA(¬+X)=1. Al tener M2(X)=1 y MAM2, por V¬ se infiere MA(¬X)=0. Por lo que, según V⊃ se puede asegurar que MA(¬+X⊃¬X)=0, lo cual significa que ¬+X⊃¬X no es válida, y según la proposición 7b, ¬+X⊃¬X∉TEO. □

Conclusión 6. (*Ubicación de LD en la jerarquía de LBVA*). El sistema LD se encuentra ubicado entre LBVA y LT en la jerarquía de lógicas basadas en LB Sierra (2010), de la siguiente manera.

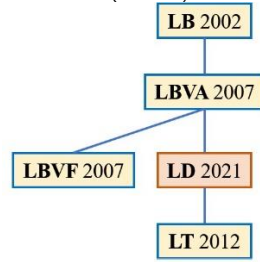

La relación: extensión estricta

Prueba. LBVA se construye con los axiomas de LC, junto a las definiciones de buena fundamentación, verdad y falsedad aristotélicas; todos ellos son teoremas de LD, por lo que LD es una extensión de LBVA. Además, según Sierra (2007a) ¬(X∪Y)⊃(¬X•¬Y) no es teorema de LBVA, y por la proposición 2b si es teorema de LD. Por lo tanto, LD es una **extensión estricta** de LBVA.

LT se construye, entre otros, con los axiomas de LC, junto con los axiomas AxR, AxT y AxMP+ de LD, y según Sierra (2012) en LT se tiene como teorema el Ax+ de LD, por lo que LT es una extensión de LD. Por otro lado, la fórmula ~+X⊃¬+X es un axioma de LT, y por el comentario 2a se sabe que ~+X⊃¬+X no es válida en LD. Por lo tanto, LT es una **extensión estricta** de LD.

La fórmula ¬+X⊃¬A, según Sierra (2007b) es válida en LBVF, pero en el comentario 2b se afirma que ¬+X⊃¬X no es válida en LD, además, la fórmula (¬X∪¬Y)⊃¬(X•Y) por la proposición 2c es teorema de LD, pero según Sierra (2007b) no es válida en LBVF. Por lo tanto, LD es **independiente** de LI. □

## 9 REFERENCIAS